\documentclass[a4paper,11pt,reqno,smallextended]{amsart} %{svjour3}

\usepackage{latexsym}
\usepackage{amssymb,amsmath,amsopn}
\usepackage{graphics}\usepackage[vcentermath]{youngtab}
\newcommand{\yngt}{\youngtabloid}

\usepackage{url}

\newtheorem{theorem}{Theorem}[section]
\newtheorem{definition}[theorem]{Definition}
\newtheorem{lemma}[theorem]{Lemma}

\newtheorem{corollary}[theorem]{Corollary}

\newtheorem{proposition}[theorem]{Proposition}

\DeclareMathOperator{\sgn}{sgn}

\DeclareMathOperator{\Hom}{Hom}

\renewcommand{\theta}{\vartheta}

\newcommand{\N}{\mathbf{N}}

\newcommand{\Ind}{\big\uparrow}

\newcommand{\ind}{\!\!\uparrow}

\newcommand{\Z}{\mathbf{Z}}

\newcommand{\oa}{1_{\raisebox{1pt}{$\scriptscriptstyle 1$}}}
\newcommand{\ob}{1_{\raisebox{1pt}{$\scriptscriptstyle 2$}}}
\newcommand{\oc}{1_{\raisebox{1pt}{$\scriptscriptstyle 3$}}}

\newcommand{\pa}{2_{\raisebox{1pt}{$\scriptscriptstyle 1$}}}
\newcommand{\pb}{2_{\raisebox{1pt}{$\scriptscriptstyle 2$}}}

\newcommand{\qa}{3_{\raisebox{1pt}{$\scriptscriptstyle 1$}}}
\newcommand{\qb}{3_{\raisebox{1pt}{$\scriptscriptstyle 2$}}}

\newcommand{\ra}{4_{\raisebox{1pt}{$\scriptscriptstyle 1$}}}
\newcommand{\rb}{4_{\raisebox{1pt}{$\scriptscriptstyle 2$}}}

\newcounter{thmlistcnt}
	{\setcounter{thmlistcnt}{0}%
	\begin{list}{\emph{(\roman{thmlistcnt})}}{%
		\usecounter{thmlistcnt}%
		\setlength{\topsep}{0pt}%
		\setlength{\leftmargin}{0pt}%
		\setlength{\itemsep}{0pt}%
		\setlength{\itemindent}{17pt}}%
	}%
	{\end{list}}%

\begin{document}

%Enlarged spacing
%\setlength{\baselineskip}{18pt}

%\journalname{Journal of Algebraic Combinatorics}

\thispagestyle{empty}

\title{Set families and Foulkes modules}

\author{Rowena Paget and Mark Wildon}

\date{March 2011}

%%\institute{Rowena Paget \at
%School of Mathematics, Statistics and Actuarial Science, 
%Cornwallis Building, University of Kent, Canterbury, 
%Kent CT2 7NF \\
% Tel.: +44 (0)1227 824755 \\
%              \email{r.e.paget@kent.ac.uk}
%           \and
%           Mark Wildon \at
%Mathematics Department, University Walk, Bristol BS8 1TW \\
%Tel.: +44 (0)117 3763820 \\
%\email{mark.wildon@bristol.ac.uk}
%}

\begin{abstract}
We construct a new family of homomorphisms from
Specht modules into Foulkes modules for the symmetric
group. These homomorphisms are used to 
give a combinatorial description
of the minimal
partitions (in the dominance order) which label the
irreducible characters appearing as summands of
the characters of Foulkes modules. The homomorphisms
are defined using certain families of subsets
of the natural numbers. 
These families are of independent
interest; we prove a number of combinatorial results
concerning them.
\end{abstract}


\maketitle

\thispagestyle{empty}

\section{Introduction}

The aim of this paper is to construct a new
family of homomorphisms from Specht modules
into Foulkes modules, and to explore some of
the consequences for the structure of Foulkes modules. 
Foulkes modules are the object of the longstanding Foulkes' Conjecture,
first made at end of \S 1
of~\cite{Foulkes},
  which spans representation theory, invariant theory and combinatorics.  We shall discuss some of these connections in this introduction.

Let~$S_r$ denote the symmetric group of degree $r$.
For $m$, $n \in \N$,
the \emph{Foulkes module}~$H^{(m^n)}$ is defined to be the 
permutation module for $\Z S_{mn}$
given by the action  of the symmetric group $S_{mn}$
on the collection of set partitions of a set of size~$mn$
into $n$ sets each of size $m$. Let $\phi^{(m^n)}$ be the permutation
character afforded by $H^{(m^n)}$. We shall refer to
these characters as \emph{Foulkes characters}.
Foulkes' Conjecture asserts that if~$m$, $n$ are natural numbers with $m < n$, 
and $\chi$
is an irreducible character of~$S_{mn}$, then the multiplicity
of $\chi$ in  $\phi^{(m^n)}$ is at least as great as 
the multiplicity of $\chi$ in $\phi^{(n^m)}$.

Foulkes' Conjecture
can be recast for $\mathrm{GL}(V)$-modules, where $V$ is a finite dimensional complex vector space. Put in these terms, it claims that if \hbox{$m < n$} 
then there is an embedding  of $\mathrm{GL}(V)$-modules $S^n(S^mV) \hookrightarrow S^m(S^nV)$, where $S^m$ denotes the $m$-th symmetric power.
In \cite[page 352]{Brion}, Brion used 
this interpretation and ideas from geometric invariant
theory to show that Foulkes conjecture is true provided $n$ is 
large compared to $m$.
%on page 352 of \cite{Brion}).

%\footnote{Howe's Conjecture~\cite[Section 2.5]{Howe} is a generalisation of Foulkes' Conjecture which considers a canonical map between the symmetric powers. M\"uller and Neunh\"offer produced a counterexample to this in~\cite{MN}.}

Finally, we mention
that taking formal characters of $\mathrm{GL}(V)$-modules gives a purely combinatorial 
formulation of Foulkes' Conjecture in terms of symmetric functions. 
In this setting it
 states that if $m < n$ then the difference of plethysms of
  Schur functions $s_n[s_m] - s_m[s_n]$ is  a sum of Schur functions with non-negative coefficients. Understanding these plethysm coefficients was identified by Stanley  in \cite[Problem 9]{Stanley} as an important
open positivity problem in algebraic combinatorics.

For a full outline of the results proved on Foulkes characters
in this article see \S 2 below.
Our main result (see Theorem~\ref{thm:minimals}) gives a combinatorial
description, in terms of certain set families, of 
the minimal partitions 
which label irreducible characters appearing as summands of Foulkes
characters. (Here minimality is with respect to the dominance order on partitions.)
This theorem gives the strongest general results on 
Foulkes characters known to date.

Using Theorem~\ref{thm:minimals}, the authors have found all 
minimal constituents of the Foulkes characters $\phi^{(m^n)}$ for
$m + n \le 20$. %{\bf [Computation is now doing $m=9$, $n=11$.]} 
The data, together
with the associated computer software,
are available from the second author's website:
\url{www.ma.rhul.ac.uk/~uvah099}. %{\bf [Will put this together.]}

It is an important feature of our approach that Theorem~\ref{thm:minimals}
is proved using an explicitly defined map between integral
Specht and Foulkes
modules for
the symmetric group. %SPECHT?
This `characteristic-free' approach is well-suited to our results,
and does not create any significant extra difficulties in their proofs.
For background on integral modules the reader is 
referred to \cite[Chapter 11]{CR}. A subsequent paper by the authors
will apply the results herein to study the behaviour of Foulkes modules
over fields of prime characteristic.

%(Here minimality is with
%respect to the dominance order on partitions, defined below.) 

Although Foulkes' Conjecture remains open, some progress has been made. 
Besides the asymptotic result of Brion already mentioned,
the conjecture is known to hold when $m \le 4$.
When Foulkes made his conjecture in 1950 it
was already known to hold by the work of Thrall (see \cite[Theorems~III
and~IV]{Thrall}) in the case $m=2$.
It was proved when $m=3$
by Dent and Siemons \cite{DentSiemons}. The most recent progress
was made in 2008 by McKay \cite{McKay1}, who proved it when $m=3$ and $m=4$.
McKay's proof uses a family of maps $\psi_{(n^m)} : H^{(n^m)}
\rightarrow H^{(m^n)}$ which were
first defined by Wagner and Siemons\footnote{See \cite{McKay1}. These maps were independently defined
by Stanley: see the discussion following Problem~9 in \cite{Stanley}. Both  Wagner and Siemons, and Stanley  formulated stronger
versions of Foulkes' Conjecture (and also of Howe's Conjecture on $\mathrm{GL}(V)$-modules; 
see \cite[Section 2.5]{Howe})
using these maps. A counterexample to these stronger conjectures 
is given in~\cite{MN}.%First counterexample in Pylyavskyy
}
in 1986.
%(McKay's paper 
%explains the connection with a later family of maps defined
%by Stanley.) 
McKay's main result is that if $\psi_{(m^m)}$ is invertible
then $\psi_{(n^m)}$ is injective for any $n > m$.
The maps $\psi_{(m^m)}$ for $m \le 4$ were already 
known to be invertible by the
work of  M{\"u}ller and 
Neunh{\"o}ffer~\cite{MN}, so Foulkes' Conjecture is proved in these cases.

The main contribution of \cite{MN}
was to show that  $\psi_{(5^5)}$ has a non-zero kernel.
 It
is however still possible that the maps $\psi_{(n^m)}$ will
have a role in proving or disproving Foulkes' Conjecture.
In particular, 
a conjecture of
Johannes Siemons\footnote{Seminar given at
Castro Urdiales, October 2007.} implies that if the kernel
of $\psi_{(n^m)}$ is non-zero, then there is a partition $\lambda$
labelling a minimal constituent of $\phi^{(n^m)}$ and a
homomorphism $f : S^\lambda \rightarrow H^{(n^m)}$ such
that $S^\lambda f \subseteq \ker \psi_{(n^m)}$.
%the image of the Specht module $S^\lambda f$ is contained in $\ker \psi_{(n^m)}$.
These minimal constituents of the Foulkes character  $\phi^{(n^m)}$ are described in our Theorem~\ref{thm:minimals}. 
It follows from
Theorem~\ref{thm:basis} that any such map $f$ is an integral
linear combination of the homomorphisms explicitly defined in \S 3.2.

%then gives a basis for the space of homomorphisms from a Specht module labelled %by a minimal constituent into the Foulkes module $H^{(n^m)}$.
%I know we haven't defined Specht modules yet but this important motivation seems to be a good thing to end the pre-intro with?

%Should integral reps comment (currently at the end of the next section) go here?

\section{Outline}

%Our results concern homomorphisms from certain Specht modules into Foulkes modules. If $\lambda$ is a partition of $r$,
%we denote by $S^\lambda$ the associated Specht module for $\Z S_r$.
%The reader is referred to \S 3 of this paper for the definition of Specht
%modules and our notational conventions for Foulkes modules.
%%(Our reasons for working with integral modules are explained below.)
%%In order to state our main results, 
Our homomorphisms are defined using certain families of 
subsets of the natural numbers. We shall need the following
combinatorial definitions.

\begin{definition}
Let $X = \{ x_1, \ldots, x_m \}$ and $Y = \{ y_1,\ldots, y_m\}$
be subsets of the natural numbers, written so that $x_1 < x_2 < \ldots
< x_m$ and $y_1 < y_2 < \ldots < y_m$. We say that $Y$ \emph{majorizes} $X$,
and write $X \preceq Y$, if $x_i \le y_i$ for each~$i$.
\end{definition}

The reader may find it helpful to refer to Figure~2 in \S 6.2
below, which shows part of the lattice of $4$-subsets of $\mathbf{N}$
under the majorization order.

\begin{definition}
A \emph{set family} of \emph{shape} 
$(m^n)$ is a 
collection of $n$ distinct $m$-subsets of the natural
numbers. A set family $\mathcal{P}$ is
\emph{closed} if whenever $Y \in \mathcal{P}$ and $X \prec Y$, then
$X \in \mathcal{P}$.
\end{definition}

In the following definition, $\lambda'$ denotes the conjugate
of the partition $\lambda$.

\begin{definition}
Let $\lambda$ be a partition with largest part of size $a$. 
A set family
has \emph{type} $\lambda$ if it has exactly $\lambda'_i$ sets
containing $i$ for each $i \in \{1,\ldots,a\}$.
\end{definition}

It is easily seen that
if $\mathcal{P}$ is a closed set family then 
 for any $i \in \N$, at least as 
many sets in $\mathcal{P}$
contain $i$ as contain $i+1$. 
Closed set families therefore
have well-defined types. The reason for working with
conjugate partitions will be seen in Theorem~\ref{thm:downset} below.

%In the following definition we use the dominance order on
%partitions. 
Recall that if $\lambda$ and $\mu$ are partitions
of $n$ then we say that $\lambda$ \emph{dominates} $\mu$,
and write $\lambda \unrhd \mu$, if $\sum_{i=1}^j \lambda_i
\ge \sum_{i=1}^j \mu_i$ for all $j \in \N$. (If $i$ exceeds
the number of parts of $\lambda$ or $\mu$, then the
corresponding part size should be taken to be $0$.)

\begin{definition}
Let $\mathcal{P}$ be a
set family of shape $(m^n)$ and type~$\lambda$.
We say that $\mathcal{P}$ is \emph{minimal} if 
there 
is no %downward 
set family $\mathcal{Q}$
of shape $(m^n)$ and type~$\mu$ with $\mu \lhd \lambda$.
\end{definition}

It is an important fact that
minimal set families are always closed;
we shall prove this fact when we first need it, 
in Proposition~\ref{prop:goingdown} in \S 5 below.

Finally, if $\lambda$ is a partition of $r$,
we denote by $S^\lambda$ the associated Specht module for $\Z S_r$.
The reader is referred to \S 3 of this paper for the definition of Specht
modules and our notational conventions for Foulkes modules.

%to the case where $m$ is odd.

We are now ready to state our main results.
Let $m$ be odd. In \S 3.2 we define  for each
closed set family $\mathcal{P}$ of shape $(m^n)$ and type $\lambda$,
a homomorphism $f_\mathcal{P} : S^{\lambda} \rightarrow H^{(m^n)}$.
%(Here $\mu'$ denotes the conjugate partition to $\mu$.)
A critical result, which we prove in \S 4 below,
is that these homomorphisms are well-defined.

\begin{theorem}\label{thm:downset}
Let $m$ be odd and let $n \in \N$. 
Let $\mathcal{P}$ be a closed
set family of shape $(m^n)$ and 
type $\lambda$. The map $f_\mathcal{P} : S^\lambda \rightarrow H^{(m^n)}$
defined in \S 3.2 is a well-defined injective
homomorphism from $S^\lambda$ to $H^{(m^n)}$.
\end{theorem}

%When defined, by change of scalars, over a field of characteristic
%zero, Specht modules are irreducible. Let $\chi^\lambda$ denote
%the irreducible character of
%the Specht module $S^\lambda$. Let $\phi^{(m^n)}$ denote 
%the character  of the Foulkes module $H^{(m^n)}$.
%It  follows 
%from Theorem~\ref{thm:downset}
%that if there is a downward closed set family of type $\lambda$
%and shape $(m^n)$, then $\chi^\lambda$ is a constituent
%of $\phi^{(m^n)}$.

Let
$\chi^\lambda$ be the
irreducible character afforded by
the Specht module $S^\lambda$. (More precisely, if~$\lambda$
is a partition of $r$, then
$\chi^\lambda$ is the character of the $\mathbf{Q}S_r$-module
$S^\lambda \otimes_{\mathbf{Z}} \mathbf{Q}$.) It is well
known that every irreducible
character of a symmetric group is equal to a $\chi^\lambda$: see, for
instance,
\cite[Theorem~4.12]{James}.
In terms of characters, Theorem~\ref{thm:downset} states
that if there is a closed set family of type $(m^n)$ and type $\lambda$
then $\left< \phi^{(m^n)}, \chi^\lambda\right> \ge 1$.

If $\pi$ is a character of $S_r$ and $\lambda$ is a partition 
of $r$, then we shall say that~$\chi^\lambda$
is a \emph{minimal constituent} of
$\pi$ if $\left<\pi, \chi^\lambda\right> \ge 1$ and
$\left<\pi, \chi^\mu\right> = 0$ if $\mu \lhd \lambda$.
In~\S 5 we prove the following theorem
which characterizes minimal constituents
of Foulkes  characters.

\begin{theorem}\label{thm:minimals}
Let $m, n \in \N$. 

\emph{(i)} If $m$ is even then the unique minimal constituent
of $\phi^{(m^n)}$ is $\chi^{(m^n)}$. % and $\left< \phi^{(m^n)}, \chi^{(m^n)} \right> = 1$.

\emph{(ii)} If $m$ is odd then 
$\chi^\lambda$
is a minimal constituent of $\phi^{(m^n)}$ if and only if 
there is a minimal 
set family of shape $(m^n)$ and type $\lambda$.
%In this event,
%the multiplicity $\left< \phi^{(m^n)}, \chi^\lambda \right>$ 
%is equal to the number of set families of shape $(m^n)$ and type $\lambda'$. 
\end{theorem}

We also show that if $m$ is even then $\left< \phi^{(m^n)}, \chi^{(m^n)} \right> = 1$, and
that part (ii) of the above theorem can be sharpened as follows.

\begin{theorem}\label{thm:basis}
Let $m$ be odd and let $n \in \N$. Suppose that $\chi^\lambda$
is a minimal constituent of $\phi^{(m^n)}$. If $\mathcal{P}_1, \ldots,
\mathcal{P}_d$ are the 
set families of shape $(m^n)$ and type~$\lambda$,
then $\mathcal{P}_1, \ldots, \mathcal{P}_d$ are closed, and
% the $\mathcal{P}_i$ are closed, and
the homomorphisms $f_{\mathcal{P}_1}, \ldots, f_{\mathcal{P}_d}$
are a $\Z$-basis for $\Hom(S^\lambda, H^{(m^n)})$.
In particular,
$\left< \phi^{(m^n)},\chi^\lambda \right> = d$.
\end{theorem}

We pause to give a small example that will illustrate these theorems. 
%the two preceding
%theorems. 
We take $m = 3$ and $n = 4$. The 
three closed set families of shape~$(3^4)$ are
\begin{align*}
&\bigl\{ \{1,2,3\}, \{1,2,4\}, \{1,2,5\}, \{1,2,6\} \bigr\}, \\
&\bigl\{ \{1,2,3\}, \{1,2,4\}, \{1,2,5\}, \{1,3,4\}  \bigr\}, \\
&\bigl\{ \{1,2,3\}, \{1,2,4\}, \{1,3,4\}, \{2,3,4\} \bigr\},
\end{align*}
of types $(6,2,2,2)$, $(5,4,2,1)$ and $(4,4,4)$ respectively.
Since these partitions are incomparable in the dominance
order, and (by Proposition~\ref{prop:goingdown}) any minimal
set family is closed,
the set families above are minimal. It now follows from
Theorems~\ref{thm:minimals} and~\ref{thm:basis} that
$\phi^{(3^4)}$ has
$\chi^{(5,2,2,2)}$, $\chi^{(5,3,2,1)}$ and $\chi^{(4,4,4)}$ as
summands, each with multiplicity~$1$. Moreover, if \hbox{$\bigl< \phi^{(3^4)},
\chi^\mu \bigr> \ge 1$} then $\mu$ dominates one of these partitions.
The presence of such larger constituents cannot be detected
by the homomorphisms $f_\mathcal{P}$. 
%In fact,
%\[ 
%\begin{split}
%\phi^{(3^4)} = \chi^{(12)} + \chi^{(10,2)} + \chi^{(9,3)}
%+ \chi^{(8,4)} +  \chi^{(8,2,2)} + \chi^{(7,4,1)}
%+ \qquad\qquad \\ \qquad\qquad
%\chi^{(7,3,2)} + \chi^{(6,6)} + \chi^{(6,4,2)} + \chi^{(6,2,2,2)} + 
%\chi^{(5,4,2,1)} + \chi^{(4,4,4)}. 
%\end{split}
%\]

%Even in this small case, it takes some effort to check this
%without the use of these theorems.

Small examples of this kind
are apt to create the false impression that every closed
set family is minimal and is the unique set family of its type. 
In~\S 6.1 we prove the following theorem
which clarifies the relationship between these properties.

\begin{theorem}\label{thm:imps}
If $\mathcal{P}$ is the 
unique set family of its type, then $\mathcal{P}$ is minimal.
If~$\mathcal{P}$ is a minimal set family then $\mathcal{P}$ is closed.
%If $\mathcal{P}$ is the unique closed set family of its type, then
%$\mathcal{P}$ is minimal. Moreover, any minimal set family is closed.
 There exist closed set families that
are not minimal, and minimal set families that are not unique
for their type.
\end{theorem}

The existence of minimal set families that are not unique for their type is
of
particular significance, 
since such families demonstrate that the multiplicity $d$ in Theorem~\ref{thm:basis}
can be strictly greater than 1.
%~\hbox{$> 1$}.

Even with the help of Theorem \ref{thm:minimals},
it appears to be a difficult matter to decide, when $m$ is odd, whether a given
partition of $mn$ labels a minimal constituent of the Foulkes
character $\phi^{(m^n)}$.
In \S 6.2 we give a construction that
gives some of these partitions. We prove
that this construction 
%gives every partition labelling
%a minimal constituent 
gives every such partition
%of $\phi^{(m^n)}$ 
if and only if $n \le 5$.

We end in \S 7 by defining the generalized Foulkes characters
$\phi^\mu$  considered in \cite{McKay1}
and showing how to obtain their minimal
constituents from the minimal constituents of the $\phi^{(m^n)}$.
This section
may be read independently of the rest of the paper.

%
%%Finally, we explain that
%%we have chosen to work throughout with integral modules because
%%this `characteristic-free' setting
%%is well-suited to our results and does not create
%%any significant extra difficulties in their proofs.
%%For background on integral modules the reader is referred to \cite[Chapter 11]{CR}. When stating results for characteristic zero representations we have used the language of characters.
%The following lemma
%collects the routine results we use to pass from integral
%homomorphisms to results about characters.
%\begin{lemma}
%Let $M$ and $N$ be $\mathbf{Z}$-free $\mathbf{Z}S_r$-modules
%and let $f : M \rightarrow N$ be a module homomorphism. Let
%
%\begin{thmlist}
%\item If $f$ is injective 
%
%extra difficulties in the proofs.
%%A subsequent paper by the authors will apply the results
%%of this paper to study
%%the behaviour of Foulkes modules over fields of prime characteristic.
%{\bf [We have some results on homomorphisms, and the analogue of
%James submodule theorem already.]}

%Using Theorem~\ref{thm:minimals},

%are distinct homomorphisms linearly independent?

\section{Specht modules and homomorphisms}
In this section we recall the definition of Specht modules
as submodules of Young permutation modules and define 
the homomorphisms $f_\mathcal{P}$.

%\subsection{Notation}
The following notation simplifies these definitions
and will be found very useful in the proofs which follow.
Given a partition $\lambda$ of~$r$ with largest part of size~$a$,
let~$A(\lambda)$ be the set consisting of the symbols
$i_j$ for $1 \le i \le a$ and $1 \le j \le \lambda'_i$.
We say that $i$ is the \emph{number} and $j$ is the \emph{index}
of the symbol $i_j$.

%\emph{Numbers and indices. $A(\lambda)$ is the 'alphabet' for
%$\lambda$. Justification is that it makes it clear at 
%a glance where an entry in a set partition sits in relation to 
%a generating polytabloid of a Specht module. Define $b_t$.}

%Let $t^\lambda$ denote the $\lambda$ tableau on the set $A(\lambda)$
%whose $i$-th column has entries from the set
%\[ \{ i_j : 1 \le j \le \lambda'_i \}, \]
%arranged so that the indices increase down each column. For example,
%\[ t^{(4,4,1)} = \young(\oa\pa\qa\ra,\ob\pb\qb\rb,\oc) \]

%Index notation.

\subsection{Specht modules}
Let $\lambda$ be a partition of $r$.
A $\lambda$-\emph{tableau} is an assignment
of the elements of $A(\lambda)$ to
the boxes of the Young diagram of $\lambda$.
Given a $\lambda$-tableau $t$, we obtain the associated 
\emph{tabloid} $\mathbf{t}$ by disregarding
the order of the elements within the rows of  $t$. For example, if
\[ t = \young(\ra\pb\qb\oa,\pa\ob\qa\rb,\oc) 
 \quad \hbox{then} \quad 
\mathbf{t} = \yngt(\ra\pb\qb\oa,\pa\ob\qa\rb,\oc) 
           = \yngt(\oa\pb\qb\ra,\ob\pa\qa\rb,\oc) 
           = \ldots \quad \hbox{etc.}\]

\smallskip
\noindent 
We may identify %the symmetric group
$S_r$ with the symmetric group on $A(\lambda)$. 
%adapted to
%the partition~$\lambda$. 
After this identification is made,
the natural permutation 
action of $S_r$ on the set of $\lambda$-tableaux
gives rise to a well-defined action of $S_r$ on the set of $\lambda$-tabloids. 
We denote the
associated permutation module for $\Z S_r$
by $M^\lambda$; it is the 
\emph{Young permutation module} corresponding
to~$\lambda$. For example, $M^{(r-1,1)}$ affords the natural 
integral representation
of~$S_r$ as $r \times r$ permutation matrices.

Given a 
$\lambda$-tableau $t$, we let $C(t)$ be the subgroup of $S_r$
consisting of those elements which fix setwise the
columns of $t$. Define $b_t \in \Z S_r$ by
%\[ b_t = \sum_{g \in C(t)} g \sgn(g). \]
\[ b_t = \sum_{\tau \in C(t)} \sgn(\tau) \tau . \]
The \emph{polytabloid} corresponding
to $t$ is the element $e_t \in M^\lambda$ defined by
\[ e_t = \mathbf{t} b_t. \] % \sum_{g \in C(t)} \textbf{t} g \sgn(g). \]
The \emph{Specht module} $S^\lambda$ is defined to be
the submodule of $M^\lambda$ spanned by the $\lambda$-polytabloids.
An easy calculation shows that if $\sigma \in S_r$ then
$(e_t)\sigma = e_{t\sigma}$, and so $S^\lambda$ is cyclic, generated by
any single polytabloid. 

It follows from Theorem 4.12 of \cite{James} that
the rational $\mathbf{Q}S_n$-modules 
$S^\lambda \otimes_\mathbf{Z} \mathbf{Q}$ for $\lambda$
a partition of $n$ are irreducible, and that they afford
all the ordinary irreducible characters of $S_n$.

 %This property is critical to the definition
%of our homomorphisms.

Let $t_\lambda$ be the $\lambda$-tableau whose $i$-th column
is $i_1, \ldots, i_{\lambda'_i}$ when read from top to bottom.
Note that the elements of
$C(t_\lambda)$ permute the indices of symbols in $A(\lambda)$ while leaving
the numbers unchanged.

%\subsection{Indexed set partitions}

\subsection{Definition of the homomorphisms $f_\mathcal{P}$}
Throughout this
section, let $m,n \in \mathbf{N}$ and let $\lambda$ be a partition of $mn$.
After identifying $S_{mn}$ with the symmetric
group on the set $A(\lambda)$, the
elements of the canonical permutation basis of $H^{(m^n)}$ are given
by the following definition.

\begin{definition}
%Let $\lambda$ be a partition of $mn$.
An \emph{indexed set partition} of shape $(m^n)$ and type $\lambda$
is a set partition of $A(\lambda)$ into $n$ sets each of size $m$.
\end{definition}

Our notation allows us to pass easily from set families to 
indexed set partitions.
%basis
%elements of $H^{(m^n)}$.

\begin{definition}
Let $\mathcal{P}$ be a set family of shape $(m^n)$ and type $\lambda$.
Order the sets making up~$\mathcal{P}$ lexicographically,
so that $\mathcal{P} = \{X_1, \ldots, X_n\}$ where $X_1 < X_2 < \cdots < X_n$.
The \emph{indexed set partition associated to}~$\mathcal{P}$
is 
the indexed set partition of type $\lambda$
obtained
by appending indices to the elements of the sets $X_1, \ldots, X_n$ 
so that the elements of $X_1$ all get the index $1$, and
an element $i \in X_r$ is given the smallest index not appended to
any $i$ appearing in~$X_1, \ldots, X_{r-1}$.
\end{definition}

For example, the indexed set partition associated to the closed set family
\[ \mathcal{Q}  = \bigl\{ \{1,2,3\}, \{1,2,4\}, \{1,2,5\}, \{1,3,4\} \} \]
of type $(5,4,2,1)$ is
\[ u = \bigl\{ 
\{1_1, 2_1, 3_1\}, \{ 1_2, 2_2, 4_1 \}, \{1_3, 2_3, 5_1\}, \{ 1_4, 3_2, 4_2\} \bigr\}
\in H^{(3^4)}. \]

Since $S^\lambda$ is generated by the polytabloid $e_{t_\lambda}$, any
homomorphism from~$S^\lambda$ is determined
by its effect on $e_{t_\lambda}$. Specifically, if $f : S^\lambda \rightarrow 
M$ is a homomorphism of $S_{mn}$-modules, 
$t$ is a $\lambda$-tableau and $\sigma \in S_{mn}$ is 
such that $t_\lambda \sigma = t$, then $e_t f = (e_{t_\lambda}f) \sigma$.

\begin{definition}
Let $\mathcal{P}$ be a set family of shape $(m^n)$ and type $\lambda$.
We define $f_\mathcal{P} : S^\lambda \rightarrow H^{(m^n)}$ by
%\begin{equation*}\label{eq:defn}
$e_{t_\lambda} f_\mathcal{P} = 
 u b_{t_\lambda}$ %\sum_{g \in C(t_\lambda)} g \sgn(g) = u b_t 
%\end{equation*}
where $u$ is the indexed set partition associated to $\mathcal{P}$.
\end{definition}

If $\mathcal{Q}$ is as above
then the homomorphism $f_{\mathcal{Q}} : S^{(5,4,2,1)} \rightarrow H^{(3^4)}$ is defined on the generator $e_{t_{(5,4,2,1)}}$ of
$S^{(5,4,2,1)}$ by
\[ e_{t_{(5,4,2,1)}} \mapsto %f_\mathcal{P} = 
\bigl\{ \{1_1, 2_1, 3_1\}, \{ 1_2, 2_2, 4_1 \}, \{1_3, 2_3, 5_1\}, \{ 1_4, 3_2, 4_2\} \bigr\} 
b_{t_{(5,4,2,1)}}.
\]

We remark that while we have, for definiteness, given an explicit
scheme for passing from set families to indexed set partitions,
a different choice will at most lead to changes of sign in the
maps $f_\mathcal{P}$. For example, if in our index appending scheme,
the lexicographic order on sets is replaced with the colexicographic
order, then the homomorphism above
would instead be defined by
\begin{align*} 
e_{t_{(5,4,2,1)}} &\mapsto 
\bigl\{ \{1_1, 2_1, 3_1\}, \{ 1_2, 2_2, 4_1 \}, \{1_3, 3_2, 4_2\}, \{ 1_4, 2_3, 5_1\} \bigr\} 
b_{t_{(5,4,2,1)}} \\
%&= -e_{t_{(5,4,2,1)}} f_\mathcal{P} 
%\end{align*}
&= -\bigl\{ \{1_1, 2_1, 3_1\}, \{ 1_2, 2_2, 4_1 \}, 
\{1_3, 2_3, 5_1\}, \{ 1_4, 3_2, 4_2\} \bigr\} 
b_{t_{(5,4,2,1)}}.
\end{align*}

\section{Proof of Theorem~\ref{thm:downset}} % down-set homomorphism theorem}

For technical reasons it will be useful to deal with the
case $m=1$ separately. 
The only closed set family of shape $(1^n)$ is 
$\mathcal{P} = \bigl\{ \{1 \}, \{2 \} ,\ldots, \{ n \} \bigr\}$,
which has type $(n)$. The homomorphism $f_\mathcal{P} : S^{(n)}
\rightarrow H^{(1^n)}$ is defined by
\[ e_{t_{(n)}} \mapsto \bigl\{ \{1_1\}, \{2_1\}, \ldots, \{n_1\} \bigr\}. \]
Since $S^{(n)}$ is the trivial $FS_{n}$-module, this map
is clearly well-defined and injective.

To show 
that the homomorphisms $f_\mathcal{P}$ are well-defined 
when $m \ge 3$,
we shall use the description of the Specht module $S^\lambda$
given by Garnir relations. The following lemma states 
a suitable form of
these relations in our numbers-and-indices notation.

\begin{lemma}\label{lemma:Garnir}
Let $U$ be a $\Z$-free $\Z S_r$-module, let $\lambda$ be a partition of $r$
and let~$t = t_\lambda$.
If $u \in U$ is such that
\[ u b_t \sum_{\sigma \in S_{X \cup Y}}  \sigma  \sgn(\sigma)= 0 \]
for every pair of subsets 
\begin{align*}
X &\subseteq \{i_j : 1 \le j \le \lambda'_i\}, \\
Y &\subseteq \{(i+1)_j : 1 \le j \le \lambda'_{i+1}\}
\end{align*} 
such that
$|X| + |Y| > \lambda'_i$, 
then there is a homomorphism
of $\Z S_r$-modules $f : S^\lambda \rightarrow U$ 
such that $e_t f = u b_t$.
\end{lemma}

\begin{proof}
%The Specht module $S^\lambda$ is generated by the polytabloid
%$e_t$. 
It follows from the remark at the top of page 102 of \cite{FultonYT} 
that 
the kernel of the surjective map $\Z S_r \rightarrow S^\lambda$
defined by $x \mapsto e_t x$ is generated, as a right $\Z S_r$-ideal,
by elements of the following two types:

%\smallskip

%\vbox{
\begin{itemize}
\item[$\bullet$] $1 - (i_j, i_k)$ for $i_j, i_k \in A(\lambda)$;
\item[$\bullet$] $G_{X,Y} = \sum \sigma \sgn(\sigma)$, where $X$ and $Y$ are as in the statement of the lemma and the sum
is over a set of right-coset representatives for the
cosets of $S_X \times S_Y$ in $S_{X \cup Y}$.
\end{itemize} %}

Clearly $u b_t$ is killed by elements of the first type, so to
prove the lemma, it will suffice to show that $u b_t G_{X,Y} = 0$
for each $G_{X,Y}$. As in the proof
of Theorem~7.2 in \cite{James}, we set 
$\overline{S}_Z = \sum_{ \sigma \in Z} \sigma \sgn(\sigma) $
for a subset $Z$ of $S_r$. 
%$\overline{S}_X$, $\overline{S}_Y$, $\overline{S}_{X \cup Y}$
%denote the signed sums of the elements in $S_X$, $S_Y$ and $S_{X \cup Y}$
%respectively. 
%As in the proof of Theorem~7.2 in \cite{James}, we let
%$\overline{S}_X$, $\overline{S}_Y$, $\overline{S}_{X \cup Y}$
%denote the signed sums of the elements in $S_X$, $S_Y$ and $S_{X \cup Y}$
%respectively. 
Note that $\overline{S}_X \overline{S}_Y G_{X,Y} =
\overline{S}_{X \cup Y}$. By hypothesis $ub_t \overline{S}_{X \cup Y} = 0$,
so we have
\[ |X|! \, |Y|! \, u b_t G_{X,Y}= 0.\]
Since $U$ is assumed to be free as a $\mathbf{Z}$-module, it follows
that $u b_t G_{X,Y} = 0$, as required.
\end{proof}

%An easy calculation shows that if $g \in S_n$ then
%$e_t g = e_{tg}$.
%Hence the homomorphism $\theta$ given by
%Lemma~\ref{lemma:Garnir} is  defined on all of $S^\lambda$ by
%$e_{tg} \theta = e_t \theta g$. 

To show that the homomorphisms $f_\mathcal{P}$
are well-defined it  suffices to check that $e_{t_\lambda} f_\mathcal{P}$ satisfies
the relations in the previous lemma. 

%This follows at
%once from the following proposition.

\begin{proposition}
Let $m \ge 3$ be odd and let $n \in \N$. 
Suppose that $\mathcal{P}$ is a closed set family of shape $(m^n)$
and type $\lambda$ and that $u$ is the indexed set partition
associated to $\mathcal{P}$. Let $t = t_\lambda$.
%Suppose 
%that $u$ is a indexed
%set partition of shape $(m^n)$ and type $\lambda$
%such that when the indices are moved from elements of sets in $u$,
%a closed set family is obtained.
If $X$ and $Y$ are as in the statement
of Lemma~\ref{lemma:Garnir} then
\[ u b_t \sum_{\sigma \in S_{X \cup Y}} \sigma \sgn(\sigma) = 0.\]
\end{proposition}

\begin{proof}
Let $\tau \in C(t)$. Suppose that
there exist $i_x \in X$ and $(i+1)_y \in Y$ 
which appear in the same set in $u\tau$. Then $u \tau (1 - (i_x, (i+1)_y)) = 0$,
and, taking coset representatives for $\left< (i_x, (i+1)_y ) \right>$
in $S_{X \cup Y}$, we see that $u\tau \sum_{\sigma \in S_{X \cup Y}}
 \sigma \sgn(\sigma)=0$.

It therefore suffices to show that 
\begin{equation}\label{eq:cancel} 
u \sum_{\tau \in C'} 
\tau \sgn(\tau) \sum_{\sigma \in S_{X \cup Y}} \sigma\sgn(\sigma)= 0
\end{equation}
where $C'$ is the set of $\tau \in C(t)$ such that no set
in $u\tau$ meets both $X$ and~$Y$. We may assume that $C'$ is non-empty.

%to which the first paragraph
%\emph{does not} apply.

Let $\theta \in C'$ and let $v= u\theta$.
None of the $|Y|$ sets in $v$ meeting $Y$ can contain an element of $X$.
%containing
%an element $(i+1)_y \in Y$ can contain an element of $X$.
At most $\lambda_i' - |X|$ of them can contain an element of the
complementary set 
$X' = \{i_{x'} : 1 \le x' \le \lambda_i', i_{x'} \not \in X \}$.
By hypothesis $\lambda_i' - |X| < |Y|$.
Hence if there are $s$ sets which meet both $Y$ and $X'$ then,
after casting out these sets, we are left with at least
$|Y| - s > |Y| + |X| - \lambda_i'$ sets which meet $Y$ but not $X'$.
Since $\mathcal{P}$ is closed,
for each such set
\[ B = \{c(1)_{b(1)}, c(2)_{b(2)}, \ldots, c(m-1)_{b(m-1)}, (i+1)_y \} \]
in $v$, there is a corresponding set
\[ A = \{c(1)_{a(1)}, c(2)_{a(2)}, \ldots, c(m-1)_{a(m-1)}, i_z \} \]
which also appears in $v$. Note that the indices $a(1), \ldots, a(m-1), z$
are determined by the numbers $c(1), \ldots, c(m-1), i$.
Since $s$ of the elements of $X'$ appear in sets which also meet $Y$,
at most $|X'| - s$ of the sets $A$ can
have $i_z \in X'$. Hence at least 
$(|Y| - s) - (|X'| - s) = |Y| - |X'|
= |Y| + |X| - \lambda_i'$
of the sets $A$ 
have $i_z \in X$. Therefore we may find sets $B$ and $A$ in $v$ so that
\begin{align*}
B &= \{c(1)_{b(1)}, c(2)_{b(2)}, \ldots, c(m-1)_{b(m-1)}, (i+1)_y \} \\
A &= \{c(1)_{a(1)}, c(2)_{a(2)}, \ldots, c(m-1)_{a(m-1)}, i_x \} 
\end{align*}
where $(i+1)_y \in Y$ and $i_x \in X$.

%Let $\theta \in C'$ and let $v= u\theta$.
%None of the $|Y|$ sets in $v$ meeting $Y$ can contain an element of $X$.
%%containing
%%an element $(i+1)_y \in Y$ can contain an element of $X$.
%At most $\lambda_i' - |X|$ of them can contain an element of the
%complementary set 
%$X' = \{i_{x'} : 1 \le x' \le \lambda_i', i_{x'} \not \in X \}$.
%By hypothesis $\lambda_i' - |X| < |Y|$. 
%Hence the number of sets in $v$ which meet $Y$
%is strictly greater than the number of sets in $v$ which meet $X'$.
%It follows, on casting out the
%sets meeting both $X'$ and $Y$, that there is a set
%\[ B = \{c(1)_{b(1)}, c(2)_{b(2)}, \ldots, c(m-1)_{b(m-1)}, (i+1)_y \} \]
%in $v$ which meets $Y$ but not $X'$. Moreover, since $\mathcal{P}$ is closed,
%we may choose~$B$ so that, 
%for some indices $a(1), \ldots, a(m-1), x$, 
%the set
%\[ A = \{c(1)_{a(1)}, c(2)_{a(2)}, \ldots, c(m-1)_{a(m-1)}, i_x \} \]
%also appears in $v$. Note that the indices $a(1), \ldots, a(m-1), x$
%are determined by the numbers $c(1), \ldots, c(m-1), i$.

Let
\[ \pi = (c(1)_{a(1)}, c(1)_{b(1)}) \cdots (c(m-1)_{a(m-1)}, c(m-1)_{b(m-1)}). \]
Since $B (i_x, (i+1)_y) = A \pi$ we have
\begin{equation}\label{eq:equal} 
v(i_x, (i+1)_y) = v\pi.
\end{equation} 
No set in $v\pi$
meets both $X$ and $Y$, so since $v\pi = u\theta \pi$, we have 
$\theta \pi \in C'$.
Thus $u \theta$ and $u\theta \pi$ are distinct summands of
$u \sum_{\tau \in C'} \tau \sgn(\tau)$, appearing with
the same sign. (This is the only
point where we use our hypotheses on $m$.)
If $\sigma_1, \ldots, \sigma_s$ is a set of right coset representatives
for the cosets of $\left< (i_x, (i+1)_y) \right>$ in $S_{X \cup Y}$ then, by \eqref{eq:equal},
\[
(u\theta + u\theta \pi) \!\!\sum_{\sigma \in S_{X \cup Y}} \sigma
\sgn (\sigma) =
v (1 + \pi) (1 - (i_x, (i+1)_y)) \sum_{r=1}^s \sigma_r 
\sgn(\sigma_r) = 0.\]

Let $H$ be the subgroup of $C(t)$ of elements
that fix all the entries in columns
$i$ and $i+1$ of $t$. 
We have shown that given any $\theta \in C'$, there exists
a non-identity
even permutation $\pi_\theta \in H$ such that $\theta \pi_\theta \in C'$
and the contributions
to~\eqref{eq:cancel} from $u \theta$ and $u\theta \pi_\theta$ cancel.

Let $C'_{i,i+1}$ be the subset of $C'$ of elements which only move entries
in columns $i$ and $i+1$ of $t$. Let $\theta, \theta' \in C'_{i,i+1}$.
%Suppose 
%the sets $\{\theta, \theta \pi_\theta\}$, 
%$\{\theta', \theta' \pi_{\theta'}\}$ meet.
%Then we must have $\theta = \theta'$, for 
If
$\theta \pi_\theta = \theta'$ then
$\theta' \theta^{-1} = \pi_\theta$; since $\pi_\theta$ fixes
the entries in columns $i$ and $i+1$ of $\theta$, this implies that
$\theta' = \theta$. Hence if the sets
 $\{\theta, \theta \pi_\theta\}$, 
$\{\theta', \theta' \pi_{\theta'}\}$ meet then $\theta = \theta'$.
We may therefore pair up the elements
of $C'_{i,i+1}$ to show that
\[ u \sum_{\theta \in C'_{i,i+1}} \theta \sgn(\theta)
 \sum_{\sigma \in S_{X \cup Y}} \sigma \sgn(\sigma) = 0.\]

%Every element
%of $C'$ can be expressed uniquely in the form $\tau \theta$
%where $\tau \in H$ and $\theta \in C'_{i,i+1}$. 
%Hence 
There exist $\tau_1, \ldots, \tau_k \in H$ such that
\[ C' = \tau_1 C'_{i,i+1} \cup \cdots \cup \tau_k C'_{i,i+1} \]
where the union is disjoint. Hence the left-hand-side of~\eqref{eq:cancel}
is
\[
\begin{split} u \sum_{r=1}^k \tau_r \sgn(\tau_r) \sum_{\theta \in C'_{i,i+1}}
\theta \sgn (\theta) \sum_{\sigma\in S_{X\cup Y}} \sigma \sgn(\sigma)
= \hskip 1in\\ \hskip 1in \Bigl( u \sum_{\theta \in C'_{i,i+1}}
\theta \sgn (\theta) \sum_{\sigma\in S_{X\cup Y}} \sigma \sgn(\sigma) \Bigr)
\sum_{r=1}^k \tau_r \sgn(\tau_r) = 0,
\end{split}
\]
as we required.
\end{proof}

To complete the proof of Theorem~\ref{thm:downset}, we must show
that 
the homomorphisms $f_\mathcal{P}$ are injective. This follows
from the following general result.
%Let
% $\lambda$
%be a partition of $r$ and let $M$ be a $\mathbf{Z}$-free 
%$\mathbf{Z}S_r$-module. If $f : S^\lambda \rightarrow M$
%is any non-zero $\mathbf{Z}S_r$-module homomorphism then $f$ induces
%a non-zero map
%\[ f' : S^\lambda \otimes_{\mathbf{Z}} \mathbf{Q}
%\rightarrow M \otimes_{\mathbf{Z}} \mathbf{Q}. \]
%Since $S^\lambda \otimes_{\mathbf{Z}} \mathbf{Q}$
%is irreducible (see \cite[Theorem~4.12]{James}),  $f'$ is injective.
%Hence the original map $f$ is also injective.
%In particular, the maps $f_\mathcal{P}$ are injective.

\begin{lemma}\label{lemma:tech}
Let $\lambda$ be a partition of $r$ and let 
$M$ be a $\mathbf{Z}$-free $\mathbf{Z}S_r$-module
with character $\pi$. If $f : S^\lambda \rightarrow M$ is
a non-zero homomorphism of $\mathbf{Z}S_r$-modules then
$f$ is injective. %and $\left< \pi, \chi^\lambda \right> \not= 0$.
\end{lemma}

\begin{proof}
The homomorphism $f$ induces a non-zero homomorphism
\[ f' : S^\lambda \otimes_{\mathbf{Z}} \mathbf{Q}
\rightarrow M \otimes_{\mathbf{Z}} \mathbf{Q}. \]
Since $S^\lambda \otimes_{\mathbf{Z}} \mathbf{Q}$
is irreducible (see \cite[Theorem~4.12]{James}),  $f'$ is injective.
Hence the original map $f$ is also injective.
%Moreover,
%the image of $f'$ is a submodule of 
%$M \otimes_{\mathbf{Z}} \mathbf{Q}$ affording
%the irreducible character $\chi^\lambda$.
%Hence 
% $\left< \pi, \chi^\lambda \right> \ge 1$.
\end{proof}

\section{Minimal constituents of Foulkes characters}

In this section we prove Theorems~\ref{thm:minimals} and~\ref{thm:basis}
on the minimal constituents of the Foulkes
characters $\phi^{(m^n)}$. 

\subsection{Even case}
Let $m$ be even. To prove part (i) of Theorem~\ref{thm:minimals},
we must show that the unique minimal constituent 
of the Foulkes character $\phi^{(m^n)}$ is $\chi^{(m^n)}$.
We do this using properties of the Young permutation modules
defined in \S 2.1. 

Let $q : M^{(m^n)}\rightarrow H^{(m^n)}$ be the map which sends 
an $(m^n)$-tabloid to the indexed set partition whose member sets
are the rows of the tabloid. 
It is clear that~$q$ is a surjective
homomorphism
of $\Z S_{mn}$-modules. 
% and that
%\[ \mathbf{t}_{(m^n)} q = \bigl\{ \{ 1_1, 2_1, \ldots, m_1 \},
% \ldots, \{ 1_m, 2_m, \ldots, n_m \}. \]
It easily follows that if~$\pi^{(m^n)}$ is the  character of
$M^{(m^n)}$ then
\[ \left< \pi^{(m^n)}, \chi^\lambda \right> \ge \left< \phi^{(m^n)}
, \chi^\lambda \right> \]
for every partition $\lambda$ of $mn$.

The constituents of the character $\pi^{(m^n)}$ are given by Young's
rule: see \cite[Chapter 14]{James}. Young's rule implies that
$\chi^{(m^n)}$ is the unique minimal constituent of $\pi^{(m^n)}$,
and that $\left< \pi^{(m^n)}, \chi^{(m^n)} \right> = 1$. 
The following lemma implies that there is a non-zero homomorphism
$S^\lambda \rightarrow H^{(m^n)}$, and hence % by Lemma~\ref{lemma:tech},
$\left< \phi^{(m^n)}, 
\chi^{(m^n)} \right> \ge 1$. Therefore $\chi^{(m^n)}$ is the 
unique minimal constituent of $\phi^{(m^n)}$ and
$\left< \phi^{(m^n)}, \chi^{(m^n)}\right> = 1$, as required.

\begin{lemma}\label{lemma:evenhom}
The Specht module $S^{(m^n)}$ is not contained in the kernel of $q$.
\end{lemma}

\begin{proof}
Let $t = t_{(m^n)}$ and let $\mathbf{t} \in M^{(m^n)}$ 
be the corresponding tabloid. By definition of $q$, 
\[ \mathbf{t} q = \bigl\{ \{ 1_1, 2_1, \ldots, m_1 \},
       \ldots, \{ 1_n, 2_n, \ldots, m_n \} \bigr\}. \]
Let $H \le C(t)$ be the subgroup of the column permuting
group of $t$ that permutes as blocks for its action the
rows of $t$. (As an abstract group, $H \cong S_n$.) For example,
if $m=4$ and $n=3$, then $H$ is generated by
$(1_1, 1_2)(2_1, 2_2)(3_1, 3_2)(4_1,4_2)$ and
$(1_1, 1_2, 1_3)(2_1, 2_2,2_3)(3_1, 3_2,3_3)(4_1,4_2,4_3)$.

Note
that since $m$ is even, every element of $H$ is an even permutation.
Let $K \le C(t)$  be the subgroup of permutations which fix the
elements $1_1, \ldots, 1_n$ in the first column of $t$.
It is easy to see that $C(t) = H K$, and so
%\[ e_t q = (\mathbf{t} q) \sum_{g \in C(t)} g \sgn(g) 
%              = (\mathbf{t} q) \sum_{h \in H} h \sum_{k \in K} k \sgn(k) 
%              = n! \, (\mathbf{t} q) \sum_{k \in K} k \sgn(k).
%\]
\[ e_t q = (\mathbf{t} q) \sum_{\tau \in C(t)} \tau \sgn(\tau) 
              = (\mathbf{t} q) \sum_{\pi \in H} \pi \sum_{\sigma \in K} 
                \sigma \sgn(\sigma) 
              = n! \, (\mathbf{t} q) \sum_{\sigma \in K} \sigma \sgn(\sigma).
\]
The summands on the right-hand-side are distinct basis elements
of $H^{(m^n)}$, hence $e_t q \not= 0$. %, as required.
\end{proof}

%It is worth noting that 
%a small extension of the argument of Lemma~\ref{lemma:evenhom} 
%shows that when $m$ is odd, $\chi^{(m^n)}$ is not a constituent
%of $\phi^{(m^n)}$. {\bf [Remark about Young filtrations?]}

\subsection{Odd case}
We start with the following
general form for a homomorphism from a Specht module
into a Foulkes module $H^{(m^n)}$ when $m$ is odd.
%It should be noted that the set families in this
%proposition are not necessarily closed.

\begin{proposition}\label{prop:hom}
Let $m$ be odd and let $n \in \N$. Let $\lambda$ be a partition
of $mn$ and let $f : S^\lambda \rightarrow H^{(m^n)}$ be a homomorphism
of $\Z S_{mn}$-modules. Let $t = t_\lambda$. There exist set families
$\mathcal{P}_1, \ldots, \mathcal{P}_k$ of shape $(m^n)$
and type $\lambda$ and integers
$a_1, \ldots, a_k$ such that 
\[ e_t f = a_1 u_1 b_t + \cdots + a_k u_t b_t \]
where $u_i$ is the indexed set partition associated to $\mathcal{P}_i$. 
\end{proposition}

\begin{proof}
Let 
\[\label{eq:coeffs} e_t f = \sum_u c_u u 
\]
where the sum is over all indexed set partitions $u$ of shape $(m^n)$
and type $\lambda$ and $c_u \in \Z$.

Let $u$ be such an indexed set partition.
If the symbols $i_j$ and $i_k$ appear in the
same set in $u$ then, since $e_t (i_j, i_k) = -e_t$
and $u (i_j, i_k) = u$, we must have $c_u = 0$.
Now suppose that $u$ contains two sets
\begin{align*} 
&\{ c(1)_{a(1)}, c(2)_{a(2)}, \ldots, c(m)_{a(m)} \},  \\
&\{ c(1)_{b(1)}, c(2)_{b(2)}, \ldots, c(m)_{b(m)} \}
\end{align*}
which become equal if the indices are removed. Let 
\[ \tau = (c(1)_{a(1)}, c(1)_{b(1)}) %(c(2)_{a(2)}, c(2)_{b(2)}) 
          \ldots
          (c(m)_{a(m)}, c(m)_{b(m)}). \]
Since $m$ is odd, $\tau$ is an odd permutation in $C(t)$.
Hence $e_t \tau = -e_t$ and $u \tau = u$, and again
we have $c_u = 0$.

These remarks show that if $c_u \not= 0$ then removing
indices from the symbols in the sets making up $u$ 
leaves a set family of shape $(m^n)$ and
type $\lambda$. If removing indices from $u$ and $v$ 
gives the same set family then $v = u \tau$ for some $\tau \in C(t)$.
Since $e_t \tau = \sgn(\tau) e_t$ we have
 $c_v = \sgn(\tau) c_u$. The proposition follows.
\end{proof}

We also need a corollary of the following combinatorial proposition,
which will be used again in \S 5 below.

%The proof of this proposition requires a further extension of
%the notion of type so that it applies to an arbitrary set family.
%We shall use this result again in \S 5.

\begin{proposition}\label{prop:goingdown}
If $\mathcal{P}$ is a minimal set family then $\mathcal{P}$ is
closed.
\end{proposition}

\begin{proof}
It will be necessary in this proof to extend the definition
of type so that it applies to all set families. We define 
the \emph{conjugate type}
of a set family~$\mathcal{P}$
of shape $(m^n)$ to be the composition $\alpha$ such that $\alpha_i$
is the number of sets in $\mathcal{P}$ containing~$i$. Note
that if $\alpha$ is a partition, then $\alpha'$ is the type (in
the usual sense)
of $\mathcal{P}$.

Let $\mathcal{P}$ be a  set family of shape $(m^n)$ and
type $\lambda$.
Suppose that $\mathcal{P}$ is not closed. We may
find $A \in \mathcal{P}$ and $i+1 \in \N$ such that $i+1 \in A$ and
the set
%\[ B = \bigl( A \cup \{i\} \bigr) \,\backslash\, \{i+1\} \]
\[ B  = A \backslash \{i+1\} \thinspace \cup \thinspace \{i\} \]
is not in $\mathcal{P}$. Let $\mathcal{Q}$ be the set family obtained
from $\mathcal{P}$ by removing $A$ and adding $B$. 
If $\alpha = \lambda'$, then the conjugate type of $\mathcal{Q}$
is $\beta$ where $\beta_i = \alpha_i + 1$, $\beta_{i+1} = \alpha_{i+1} - 1$
and $\beta_j = \alpha_j $ if $j\not= i,i+1$. Hence $\beta \rhd \alpha$
(where $\rhd$ now refers to the dominance order on compositions).
Iterating this construction, we will reach a closed
set family $\mathcal{R}$ of conjugate type $\gamma$ where $\gamma \rhd \alpha$.
Since~$\mathcal{R}$ is closed, $\gamma$ is a partition. If $\nu = \gamma'$
then, $\nu' \rhd \lambda'$, and so $\nu \lhd \lambda$. Thus $\mathcal{R}$
has smaller type than $\mathcal{P}$ and so $\mathcal{P}$ is not minimal.
\end{proof}

\begin{corollary}\label{cor:tech}
If $\left<\phi^{(m^n)}, \chi^\lambda \right> \ge 1$ then 
there is a minimal %closed 
set family of shape $(m^n)$ and type $\mu$
where $\mu \unlhd \lambda$, and this set family is closed.
\end{corollary}

\begin{proof}
The hypothesis implies that there is an injective
homomorphism of $\mathbf{Q}S_{mn}$-modules 
\[ f : S^\lambda \otimes_\mathbf{Z}
\mathbf{Q} \rightarrow H^{(m^n)} \otimes_\mathbf{Z} \mathbf{Q}. \] 
Let
\[ e_{t_\lambda} f = \sum_u c_u u \]
where $c_u \in \mathbf{Q}$ and the sum is over all indexed
set partitions $u$ of shape $(m^n)$ and type $\lambda$. For each such
$u$, let $c_u = a_u/b_u$ where $a_u, b_u \in \mathbf{Z}$. Let 
$m$ be the product
of all the $b_u$. It is easy to see that the map $g : S^\lambda \rightarrow H^{(m^n)}$
 defined by
\[ e_{t_\lambda} g =  \sum_u mc_u u \]
is a well-defined injective homomorphism of $\mathbf{Z}S_{mn}$-modules. 

Applying Proposition~\ref{prop:hom} to $g$, we see
that there 
is a set family of shape $(m^n)$ and type~$\lambda$. 
Proposition~\ref{prop:goingdown} implies that
a set family of minimal type $\lhd\, \lambda$ is closed.
\end{proof}

%\begin{proof}
%By hypothesis there is a non-zero homomorphism of $\Z S_{mn}$
%modules $S^\lambda \rightarrow H^{(m^n)}$. It follows from 
%Proposition~\ref{prop:hom} that there is a set family $\mathcal{P}$
%of shape $(m^n)$ and type $\lambda$. Hence, by
%Proposition~\ref{prop:goingdown}, there is a closed set family
%of type $\nu$ where $\nu \unlhd \lambda$. The corollary follows
%by taking a minimal closed set family of type $\mu \unlhd \nu$.
%\end{proof}

We are now ready to prove part (ii) of Theorem~\ref{thm:minimals}.
Suppose
that $\chi^\lambda$ is a minimal constituent of $\phi^{(m^n)}$. 
By Corollary~\ref{cor:tech} there is a minimal set family~$\mathcal{Q}$
of type $\mu$ where $\mu \unlhd \lambda$. Since $\mathcal{Q}$
is closed, it follows from Theorem~\ref{thm:downset} that
there is a non-zero homomorphism $S^\mu \rightarrow H^{(m^n)}$, and so
%by Lemma~\ref{lemma:tech},
$\left< \phi^{(m^n)}, \chi^\mu \right> \ge 1$. Therefore
$\mu = \lambda$ and $\mathcal{Q}$ is a minimal set family of
type $\lambda$.

Conversely, suppose that there is a minimal
set family $\mathcal{P}$
of type~$\lambda$. By Proposition~\ref{prop:goingdown}, $\mathcal{P}$
is closed, and so it follows from Theorem~\ref{thm:downset}
that $\chi^\lambda$ is a summand
of $\phi^{(m^n)}$. If there is a summand $\chi^\mu$ of 
$\phi^{(m^n)}$ with $\mu \lhd \lambda$, then
by Proposition~\ref{prop:hom}, there is a set family of 
type $\mu$; this
contradicts the minimality of $\mathcal{P}$. Therefore $\chi^\lambda$
is a minimal constituent of $\phi^{(m^n)}$. 
%This completes
%the proof of part (ii) of Theorem~\ref{thm:minimals}.

\subsection{Proof of Theorem~\ref{thm:basis}}

%Suppose that $\chi^\lambda$ is a minimal constituent of $\phi^{(m^n)}$.
%By Theorem~\ref{thm:minimals}, $\lambda$ is the minimal type of
%a closed set partition of shape $(m^n)$. Hence, by
%Proposition~\ref{prop:goingdown}, if $\mathcal{Q}$ is a 
%set family of shape $(m^n)$ and type $\lambda$ then $\mathcal{Q}$ 
%is closed. 
%It now follows from Proposition~\ref{prop:hom}
%that if $f : S^\lambda \rightarrow H^{(m^n)}$ is a homomorphism of 
%$\Z S_{mn}$-modules, then there exist $a_1, \ldots, a_d \in \Z$
%such that
%\[ f = a_1 f_{\mathcal{P}_1} + \cdots + a_d f_{\mathcal{P}_d} \]
%where $\mathcal{P}_1, \ldots, \mathcal{P}_d$ are the set families
%of shape $(m^n)$ and type $\lambda$.

Suppose that $\chi^\lambda$ is a minimal constituent of~$\phi^{(m^n)}$.
By Theorem~\ref{thm:minimals}, $\lambda$ is the type of a minimal
set family of shape~$(m^n)$.  Let $u_1, \ldots, u_d$ be the indexed
set partitions associated to the
set families $\mathcal{P}_1, \ldots, \mathcal{P}_d$ of shape $(m^n)$
and type $\lambda$. Proposition~\ref{prop:goingdown} implies
that the $\mathcal{P}_r$
are closed.
By Proposition~\ref{prop:hom}
we know that if $f : S^\lambda \rightarrow H^{(m^n)}$ is a homomorphism of 
$\Z S_{mn}$-modules, then there exist $a_1, \ldots, a_d \in \Z$
such that
\[ e_t f = a_1 u_1 b_t + \cdots + a_d u_d b_t \]
where $t = t_\lambda$.
Hence %Since the $\mathcal{P}_i$ are closed, we have
\[ f = a_1 f_{\mathcal{P}_1} + \cdots + a_d f_{\mathcal{P}_d}. \]

To show that  homomorphisms $f_{\mathcal{P}_1}, \ldots, f_{\mathcal{P}_d}$
are linearly independent it suffices to show that the images of $e_t$,
\[ e_t f_{\mathcal{P}_r} = \sum_{\tau \in C(t)} u_r \tau \sgn(\tau) \]% \in H^{(m^n)} \]
are linearly independent. Given $\tau \in C(t)$ and $r \in \{1,\ldots, d\}$,
we can easily recover $\mathcal{P}_r$ from $u_r \tau$
by removing the indices from the
symbols in the sets making up $u_r \tau$. Therefore each 
$e_t f_{\mathcal{P}_r}$ is a sum of different
basis elements of $H^{(m^n)}$;
as such, they are linearly independent.

\section{Set families and partitions}

%In this section we prove some results with a more combinatorial
%flavour on set families and partitions.

\subsection{Minimal and unique set families}
We proved in Proposition~\ref{prop:goingdown} that minimal
set families are closed. The following proposition
implies that if $\mathcal{P}$ is the unique
set family of its shape and type, then $\mathcal{P}$
is minimal. This establishes the following 
chain of implications
on set families of a given shape:
\begin{equation}
\label{eq:imps}
\text{unique of its type} \implies \text{minimal} \implies \text{closed}. 
\end{equation}

\begin{proposition}
If $\mathcal{P}$ is a set family of shape $(m^n)$ and type $\lambda$ and 
$\mu \rhd \lambda$ then there are two distinct set families
of shape $(m^n)$ and type $\mu$.
\end{proposition}

\begin{proof}
We may assume that $\lambda$ and $\mu$ are neighbours
in the dominance order, so
$\mu$ is obtained from $\lambda$ by moving a box
upwards in its Young diagram. Suppose the box is moved
from column $i$ to column $j > i$. We have
$\mu'_i = \lambda'_i - 1$, $\mu'_j = \lambda'_j + 1$ and
$\mu'_k = \lambda'_k$ if $k\not= i,j$. 

The sets in $\mathcal{P}$ either contain both $i$ and $j$,
or $i$ alone, or $j$ alone, or neither. Since $\lambda'_i - \lambda'_j
= \mu'_i - \mu'_j + 2$, at least two more sets contain $i$ alone
then contain $j$ alone.
Hence there are two sets
\begin{align*}
A &= \{ x(1), x(2), \ldots, x(m-1), i \} \in \mathcal{P}, \\ 
A' &= \{ x'(1), x'(2), \ldots, x'(m-1), i \}
 \in \mathcal{P} 
\end{align*}
such that
\begin{align*}
B &= \{ x(1), x(2), \ldots, x(m-1), j \}\not\in \mathcal{P}, \\
B' &= \{ x'(1), x'(2), \ldots, x'(m-1), j \}
\not\in \mathcal{P}.
\end{align*}
Let $\mathcal{Q}$ be the
set family obtained from $\mathcal{P}$ by removing 
$A$ and adding~$B$, and let $\mathcal{Q}'$ be the set family
obtained from $\mathcal{P}$ by removing $A'$ and adding $B'$.
Then $\mathcal{Q}$ and $\mathcal{Q}'$ are two
different set families of type $\mu$.
%$\mathcal{Q} = \mathcal{P} \backslash \{ A\} \cup\{ B \}$
%and $\mathcal{Q}' = \mathcal{P}' \backslash \{A' \cup B'$ are
%different and both have type $\mu$.
\end{proof}

%We mention without proof that a set family consisting of $2$-subsets
%of the natural numbers is closed if and only if it is unique
%for its type. 

To complete the proof of Theorem~\ref{thm:imps}, we must show that
neither of the
 implications in the chain~\eqref{eq:imps} is
reversible. We shall use the following definition
throughout the remainder of this section.
%examples.

\begin{definition}
Let $A$ be a subset of the natural numbers. The \emph{downset}
of $A$, denoted $A^{\preceq}$, is the set family consisting
of all subsets $X$ such that $X \preceq A$.
\end{definition}

It is obvious that a set family is closed if and only if
it is a union of downsets.

Our first example concerns the downset $\mathcal{P} = 
\{2,4,6,8\}^\preceq$.
This is a closed set family of shape $(4^{42})$, 
but it is not unique for its
type, since the set family
\[ Q = \mathcal{P} \backslash \bigl\{ \{2,4,6,8\}, \{1,3,5,7\} \bigr\}
                   \cup       \bigl\{ \{1,2,7,8\}, \{3,4,5,6\} \bigr\} \]
has the same shape and type  as $\mathcal{P}$. 
Neither is $\mathcal{P}$ minimal. This is most easily
seen by noting that
since $\mathcal{Q}$ contains $\{1,3,5,8\}$
but not $\{1,3,5,7\}$, $\mathcal{Q}$ is not closed. 
%Hence, by 
%Proposition~\ref{prop:goingdown}, 
%there is a closed 
%set family of smaller type than~$\mathcal{P}$ and~$\mathcal{Q}$. 
Following the proof of Proposition~\ref{prop:goingdown}
leads one to (amongst others) the closed set family
\begin{align*} 
R &= \mathcal{Q} \backslash \bigl\{ \{2,4,5,8\} \bigr\} \cup
\{\{1,3,5,7\} \bigr\} \\
 &= \mathcal{P} \backslash \bigl\{ \{2,4,6,8\}, \{2,4,5,8\} \big\}
\cup \big\{ \{1,2,7,8\}, \{3,4,5,6\} \bigr\} \\
&= \{1,2,7,8\}^\preceq \cup \{1,4,6,8\}^\preceq \cup\{2,3,6,8\}^\preceq \cup\{2,4,6,7\}^\preceq \cup\{3,4,5,6\}\}^\preceq
\end{align*}
which has smaller type than $\mathcal{P}$ and $\mathcal{Q}$.
(In fact, $\mathcal{R}$ is minimal.)

%minimal set family $\mathcal{R}$ obtained from $\mathcal{P}$ by
%removing $\{2,4,6,8\}$ and $\{1,3,5,7\}$ and adding
%$\{3,4,5,6\}$ and $\{1,3,5,7\}$. Written as a union of downsets,
%\[ \mathcal{R} = \{1,4,6,8\}^\preceq \cup 
%                 \{2,3,6,8\}^\preceq \cup
%                 \{2,4,5,8\}^\preceq \cup
%                 \{2,4,6,7\}^\preceq \cup
%                 \{3,4,5,6\}^\preceq . \]
%%(One can show that $\mathcal{R} is 

The following example of a minimal set family that is not unique
for its type was
found by a  computer search for minimal set families with prescribed
shape and maximum entry.
(A description of the algorithm used 
and accompanying source code
is available from the second author's website, \url{www.maths.bris.ac.uk/~mazmjw}.)
Let
\begin{align*}
\mathcal{S} &=  \{1,5,9\}^\preceq \cup \{1,6,8\}^\preceq \cup 
                  \{2,6,7\}^\preceq \cup
                  \{3,4,8\}^\preceq \cup \{3,5,6\}^\preceq, \\
  \mathcal{S'} &= \{1,4,9\}^\preceq \cup 
                  \{1,7,8\}^\preceq \cup 
                   \{2,3,9\}^\preceq \cup \{2,4,8\}^\preceq  \hskip 1in \\
                &   \hskip 2in
                \cup   \{2,5,7\}^\preceq \cup 
                 \{3,4,7\}^\preceq \cup \{4,5,6\}^\preceq. 
\end{align*}
The set families $\mathcal{S}$ and $\mathcal{S}'$ both have type  $\lambda = (24,19,17,16,13,12,10,8,4)'$
and shape $(3^{41})$. Using the computer to enumerate all closed
set families of shape $(3^{41})$ with maximum entry $\le 9$
confirms that there are no set families of shape $(3^{41})$
with type $\lhd\; \lambda$.

\subsection{Constructing types of minimal set families} 

Let $m,n \in \mathbf{N}$.
Given a partition $\nu$ of $n-1$ with $k \le m$ parts,
we let
$(m^n) \star \nu$ denote the partition obtained from $(m^n)$ by deleting
$\nu_i$ boxes from column $m+1-i$
and then adding $\nu_i$ boxes to row~$i$, for 
each $i$ such that $1 \le i \le k$. This
construction is illustrated in Figure~1. 
%(We
%leave it to the reader to verify that the partition
%$(m^n) \star \nu$ is well-defined.)
(If
$\nu_1 = a$ then we add boxes to rows $1, 2, \ldots, k$,
and remove boxes from rows $n-a+1, \ldots, n$; since
$k + a \le n$, the partition $(m^n) \star \nu$
is well-defined.)

%that if $\nu_1 = a$, then the
% highest row from which we remove boxes is
%row $n - a$; since $k + a \le \n-1$ we have $n-a > k$;
%the construction is therefore well-defined.

\begin{figure}[t]
\begin{center}
\qquad\qquad
\includegraphics{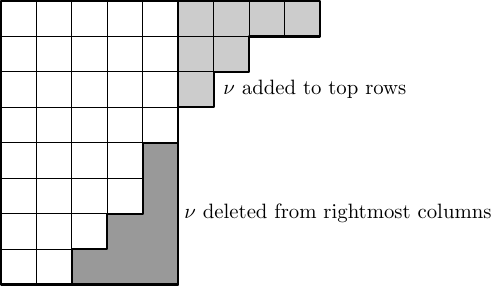}
\caption{Construction of the partition
$(5^8) \star (4,2,1)$.}
\end{center}
\end{figure}

%We shall prove that if $\nu$ is a partition of $n-1$ with
%at most $m$ parts, then there is a unique set family of 
%shape $(m^n)$ and type $(m^n) \star \nu$. 

Our first object in this section is to prove the following proposition.

\begin{proposition}\label{prop:cons} Let $\nu$ be a partition of $n-1$ with
at most $m$ parts. There is a unique set family of shape $(m^n)$
and type $(m^n) \star \nu$.
\end{proposition}
\begin{proof}
Suppose that $\nu$ has exactly $k$ parts.
Let $A = \{1,2,\ldots, m\}$.
It is easily checked that 
\[ \mathcal{P} = \bigl\{ A \bigr\} \cup 
\bigl\{ A \backslash \{ m-i + 1\} \cup \{ m + j\} :
1 \le i \le k, 1 \le j \le \nu_i \bigr\} \] 
is a set family of shape $(m^n)$ and type $(m^n) \star \nu$.

 Now suppose 
that $\mathcal{Q}$ is a set family of this shape and type.
We may write
\[ \mathcal{Q} = \bigl\{ A \backslash B_r \cup C_r : 1 \le r \le n \bigr\} \]
for some subsets $B_r,C_r \subseteq \mathbf{N}$ such that
$B_r \cap C_r = \varnothing$ for $1 \le r \le n$. 
%Observe
%that 
For each~$i$ such that $1 \le i \le k$,
exactly $\nu_i$ of the sets $B_1, \ldots, B_n$ contain $m-i+1$. 
It follows that $|B_1| + \cdots + |B_n| = n-1$, and so one
of the sets $B_r$ is empty, and the remaining $n-1$ are singletons.
Hence $A \in \mathcal{Q}$, and for each $i$ such that $1\le i\le k$, 
there are exactly $\nu_i$ sets in $\mathcal{Q}$ of the form
$A \backslash \{ m-i+1\} \cup \{ c \}$ where $c \in \{m+1, \ldots, m+\nu_1\}$.
Looking first at the case $i=1$, we see that 
$A \backslash \{ m \} \cup \{ m+ j \} \in \mathcal{Q}$ for
each $j$ such that $1 \le j \le \nu_1$. Iterating this argument for
$i = 2, \ldots, k$
shows that $\mathcal{Q} = \mathcal{P}$.
\end{proof}

We isolate the following corollary of Proposition~\ref{prop:cons}

\begin{corollary}\label{cor:mins}
Let  $m$ be odd and let $n \in \mathbf{N}$. If $\nu$
is a partition of $n-1$ with at most $m$ parts then
\[ \left< \phi^{(m^n)}, \chi^{(m^n) \star \nu} \right> =  1 \]
Moreover $\left< \phi^{(m^n)}, \chi^\mu \right> = 0$ if 
$\mu \lhd \chi^{(m^n) \star \nu}$.
\end{corollary}

\begin{proof}
By Theorem~\ref{thm:imps}, $(m^n) \star \nu$ is the type of a 
minimal set family of shape~$(m^n)$. The result now follows from 
Theorem~\ref{thm:basis}.
\end{proof}

It is natural to ask when every minimal constituent $\phi^{(m^n)}$
arises from this construction. We shall show in Proposition~\ref{prop:summary}
below that this is
the case if and only if $n \le 5$ or $m=1$. We begin with the following
straightforward lemma.

\begin{lemma}\label{lemma:prelim}
Let $\mathcal{P}$ be a closed set family of shape $(m^n)$
where $m \ge n$.
If $X \in \mathcal{P}$ then
$X \supseteq \{ 1,2, \ldots, m-n+1 \}$.
\end{lemma}

\begin{proof}
%Let $X \in \mathcal{P}$.
Suppose that the smallest number not present in $X$ is $m-t+1$, so
$X = \{1,2,\ldots, m-t, x(1), \ldots, x(t)\}$ for some $x(r) > m-t+1$.
For each $r$ such that $1 \le r \le t$, the set $X \backslash \{ x(r) \}
\cup \{ m-t+1 \}$ is majorized by $X$, so must lie in $\mathcal{P}$. Hence 
$|\mathcal{P}| \ge t+1$, and so $t \le n-1$. It follows that
$m-t \ge m-n+1$, as required.
\end{proof}

In the proof of the following
lemma, a further construction
on partitions will be found useful: 
given a partition $\lambda$ with exactly $k$ parts, each of 
size~$\ge c$, let
$\lambda - (c^k)$ denote the partition obtained
from $(\lambda_1 - c,
\ldots, \lambda_k - c)$ by removing any final parts
of size zero.

\begin{lemma}\label{lemma:all}
Let $m \in \N$ and let $n \le 5$. If $\mathcal{P}$ is a closed
set family of shape $(m^n)$ then $\mathcal{P}$ has type
$(m^n) \star \nu$ for some partition $\nu$ of $n-1$ with
at most~$m$ parts.
\end{lemma}

\begin{proof}
If $m < n$ then the result can be checked directly.
For example, when $m=4$ and $n=5$, the set family $\mathcal{P}$
must consist of $5$ sets taken from the first $5$ levels
of the lattice of $4$-subsets of $\mathbf{N}$ under the majorization
order. It is easily seen from Figure~2 that there are 
$5$ possibilities for $\mathcal{P}$, namely
$\{1,2,3,8\}^{\preceq}$, $\{1,2,3,7\}^{\preceq} \cup
\bigl\{ \{1,2,4,5\} \bigr\}$,
$\{1,2,4,6\}^{\preceq}$,
$\{1,3,4,5\}^\preceq \cup \bigl\{ \{1,2,3,6\} \bigr\}$ and 
$\{2,3,4,5\}^\preceq$.
%\[ \{1,2,3,8\}^{\preceq}, \{1,2,3,7\}^{\preceq} \cup
%\bigl\{ \{1,2,4,5\} \bigr\},
%\{1,2,4,6\}^{\preceq},
%\{1,3,4,5\}^\preceq \cup \bigl\{ \{1,2,3,6\} \bigr\} \quad\text{and}\quad 
%\{2,3,4,5\}^\preceq.
%\]
The types of these set families
are $(4^5) \star \nu$ where $\nu = (4)$, $(3,1)$,
$(2,2)$, $(2,1,1)$ and $(1^4)$ respectively.

Now suppose that $m \ge n$.
By Lemma~\ref{lemma:prelim}, we know that every set
in $\mathcal{P}$ contains $\{1,2, \ldots, m-n+1\}$. 
Let $\mathcal{Q}$ be the set family obtained by removing the elements $1$,
$2$, \ldots, $m-n+1$
%$\{1,2, \ldots, m-n+1\}$
from every set in $\mathcal{P}$, and then subtracting
$m-n+1$ from each remaining element. The shape of $\mathcal{Q}$
is $((n-1)^n)$. If $\mathcal{P}$ has type $\lambda$, where $\lambda$
has exactly $k$ parts,
%= (\lambda_1, \ldots, \lambda_k)$ 
then
$\mathcal{Q}$ has type $\lambda - ((m-n+1)^k)$.
By the result already proved, we know that
\[ \lambda - ((m-n+1)^k) = ((n-1)^n) \star \nu \]
for some partition $\nu$ of $n-1$. Since 
$\lambda - ((m-n+1)^k)$ is a partition of $(n-1)k$,
we must have $k = n$. It is now easy to see that $\lambda = (m^n) \star \nu$,
as required.
\end{proof}

\begin{figure}[h]
\begin{center}
\includegraphics{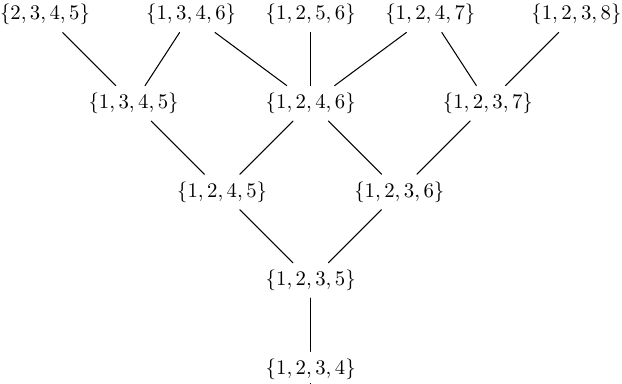}
\end{center}
\caption{The first five levels of the lattice of $4$-subsets
of~$\mathbf{N}$ under the majorization order.}
\end{figure}

%It is worth noting that Proposition~\ref{prop:cons} and Lemma~\ref{lemma:all}
%together imply that every closed set family
%of shape $(m^n)$ for $n \le 5$ is minimal and unique for its type.

We now show that Lemma~\ref{lemma:all} is false when $n \ge 6$ and $m \ge 2$.

\vbox{
\begin{lemma}\label{lemma:notmore}
Let $n \ge 6$ and let $m \ge 2$. The unique
set family of shape $(m^n)$ and type $(m+n-4, m+2, m+2, (m-1)^{n-6}, (m-2)^3)$
is 
\[ \{1,2, \ldots, m-2, m+1, m+2 \}^\preceq \cup 
   \{1,2, \ldots, m-2,m-1, m+n-4 \}^\preceq. \]
%(If $m=2$ then the second downset should be deleted.)
This type is not of the form $(m^n) \star \nu$ for any partition $\nu$
of $n-1$ with at most $m$ parts. 
\end{lemma}}

\begin{proof}
It is routine to check that this set family does have the 
claimed type. Suppose that $\mathcal{P}$ is a set family of
this type. Every set in $\mathcal{P}$ contains $\{1,2,\ldots, m-2\}$,
so, as in the proof of Lemma~\ref{lemma:all}, we may pass
to a set family $\mathcal{Q}$ of shape $(2^n)$ and type $(n-2, 4, 4, 1^{n-6})$
by removing the elements $1,2, \ldots, m-2$ from every set, and
then subtracting $m-2$ from every remaining element.
It suffices to show that $\mathcal{Q} = 
\{3,4\}^\preceq \cup \{1,n-2\}^\preceq$. Of the $n-3$ sets in~$\mathcal{Q}$
which contain $1$, at most one can contain $2$, so $2$ must appear
in two of the three remaining sets in $\mathcal{Q}$. 
A similar argument with $3$ and $4$
shows that $\{2,3\}, \{2,4\}, \{3,4\} \in \mathcal{Q}$. 
It is now clear that $\mathcal{Q}$ is as claimed.

If the final assertion in the lemma is false, then 
the type of $Q$ is  $(2^n) \star \nu$ where $\nu$ is a partition
of $n-1$ with at most $2$ parts. However, $(n-2,4,4,1^{n-6})$ has
three parts of size $> 2$, so this is impossible.
\end{proof}
 
We now use Lemmas~\ref{lemma:all} and~\ref{lemma:notmore}
to prove the following proposition.
%The following proposition follows from
%Lemmas~\ref{lemma:all} and~\ref{lemma:notmore} using 
%Theorem~\ref{thm:minimals}(ii) and Theorem~\ref{thm:imps} (which
%is needed to
%show that the unique set family of Lemma~\ref{lemma:notmore} is
%minimal).
%It requires no further proof.

\begin{proposition}\label{prop:summary}
Let $m$ be odd and let $n \in \N$. 
Every minimal constituent
of $\phi^{(m^n)}$ is of the form $\chi^{(m^n) \star \nu}$ for
some partition $\nu$ of $n-1$ with at most $m$ parts if and only
if $m=1$ or $n \le 5$.\hfill$\Box$
\end{proposition}

\begin{proof}
Since $\phi^{(1^n)} = \chi^{(n)}$ and $(n) = (1^n) \star (n-1)$,
the proposition holds when $m=1$. When $n \le 5$ it follows from
Lemma~\ref{lemma:all}, Theorem~\ref{thm:minimals} and Theorem~\ref{thm:imps}. 
If $n \ge 6$ then, since the unique
set family constructed in 
Lemma~\ref{lemma:notmore} is minimal by Theorem~\ref{thm:imps},
there is a minimal constituent of $\phi^{(m^n)}$ not of the
form $\chi^{(m^n) \star \nu}$ for any partition $\nu$ of $n-1$.
\end{proof}

\section{Minimal constituents of generalized Foulkes characters}

We end by showing how to construct the minimal constituents
of a wider class of permutation characters.
Let $\mu$ be a partition of $N$ with largest part of size $a$.
If $\mu$ has exactly 
$n(i)$ parts of length $i$ for each $i$ such that
$1 \le i \le a$, 
we define the \emph{generalized Foulkes character} $\phi^\mu$
to be the induced character
\[ \phi^{\mu}  = \Bigl( \phi^{(1^{n(1)})} \times
\phi^{(2^{n(2)})} \times \cdots \times
\phi^{(a^{n(a)})} \Bigr) \Ind_{S_{n(1)} \times S_{2n(2)}
\times \cdots \times S_{an(a)}}^{S_N} . \]
If $\mu = (m^n)$ for some $m,n \in \mathbf{N}$ then this
definition agrees with the one given earlier. 
%(For an
%equivalent definition, see \cite[\S 1]{McKay1}.)

Our main aim in this section is to prove Proposition~\ref{prop:gen}
below describing the minimal constituents of generalized Foulkes
characters. To state this result 
we need one final construction on partitions.

\begin{definition}
Given $\lambda$ a partition of $r$ and $\mu$  a partition
of $s$, we denote by $\lambda \cup \mu$ the partition of $r+s$
whose multiset of parts is the union of the multisets of parts
of $\lambda$ and $\mu$.
%Let $\lambda$ be a partition of $r$ with largest part of size $a$
%and let $\mu$ be a partition of $s$ with largest part of size $b$.
%Suppose that $\lambda$ has exactly $\alpha_i$ parts of size $i$,
%and $\mu$ has exactly $\beta_j$ parts of size $j$, for each
%$j$ such that $1 \le j < \mathrm{max}(a,b)$.
%We denote by $\lambda \cup \mu$ the partition of $r+s$ which has
%exactly $\alpha_i + \beta_i$ parts of size $i$ 
%for each
%$j$ such that $1 \le j < \mathrm{max}(a,b)$.
\end{definition}

For example, $(5,2) \cup (3,2,2,1,1) = (5,3,2,2,2,1,1)$.

\begin{proposition}\label{prop:gen}
Let $\mu = (m(1)^{n(1)}) \cup \cdots \cup (m(t)^{n(t)})$
where the $m(k)$ are distinct. 
If $\chi^\lambda$ is a minimal constituent of $\phi^\mu$
then 
\[ \lambda = \nu(1) \cup \nu(2) \cup \cdots \cup \nu(t) \]
where $\nu(k)$ is a partition of $m(k)n(k)$ and
$\chi^{\nu(k)}$ is a minimal constituent of $\phi^{(m(k)^{n(k)})}$
for each~$k$ such that $1 \le k \le t$.
\end{proposition}

Our proof of Proposition~\ref{prop:gen} uses the following
two general lemmas. % that are of wider interest.

\begin{lemma}\label{lemma:lr}
Let $\lambda$ be a partition of $r$ and let $\mu$ be a partition of 
$s$. The unique minimal constituent of $(\chi^\lambda \times \chi^\mu)
\ind_{S_r \times S_s}^{S_{r+s}}$ is $\chi^{\lambda \cup \mu}$.
\end{lemma}

\begin{proof}
Let $\theta = (\chi^\lambda \times \chi^\mu)
\ind_{S_r \times S_s}^{S_{r+s}}$. 
It follows from the description of the
Littlewood--Richardson rule given
in \cite[Chapter~16]{James} that $\chi^{\lambda \cup \mu}$ is
a constituent of the character $\theta$. A typical example, which
shows how the parts of $\lambda \cup \mu$ may be obtained 
by adding numbers to~$\lambda$, is given
in
Figure~3. 
Note that at step $j$,
the lowest $\mu_j$ positions that are eligible to be filled
receive a $j$.
 For an explanation of the notation and method used, the
reader is referred to \cite[Chapter~16]{James}.
 
%A typical example is
%shown in Figure~3 below.

The remainder of the proof can be completed using the easier
Young's rule.
The
character
\[ \psi = \bigl( \chi^\lambda \times (1_{S_{\mu}}\ind^{S_{s}}) \bigr)
\Ind_{S_r \times S_{s}}^{S_{r+s}} \]
certainly contains all the constituents of $\theta$, so to prove
the lemma, it suffices to show that $\psi$ 
has $\chi^{\lambda \cup \mu}$ as its least constituent. This
follows by induction on the number of parts of $\mu$ 
if we rewrite~$\psi$ as
\[ \bigl( \chi^\lambda \times 1_{\mu_1} \times \cdots \times 1_{\mu_k} \bigr)
\Ind_{S_r \times S_{\mu_1} \times \cdots \times S_{\mu_k}}^{S_{r + s}}\]
and then repeatedly apply Young's rule (see \cite[Chapter 14]{James}).
\end{proof}

\begin{figure}[t]
\begin{center}
\includegraphics{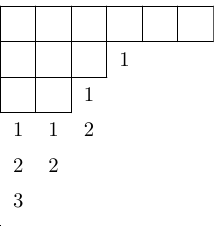}
\end{center}
\caption{
Use of the Littlewood--Richardson rule to show that
$\bigl(\chi^{(6,3,2)} \times \chi^{(4,3,1)}\bigr)
\Ind_{S_{10} \times S_8}^{S_{18}}$ has $\chi^{(6,3,2) \cup (4,3,1)}
= \chi^{(6,4,3,3,2,1)}$
as a summand.}
\end{figure}

\begin{lemma}\label{lemma:join}
Let $\pi$ be a character of $S_r$ and let $\theta$ be a character of
$S_s$. If $\chi^\nu$ is a minimal constituent of $(\pi \times \theta)
\ind_{S_r \times S_s}^{S_{r+s}}$
then $\nu = \lambda \cup \mu$ where $\chi^\lambda$ is a minimal constituent
of $\pi$ and $\chi^\mu$ is a minimal constituent of $\theta$.
\end{lemma}

\begin{proof}
Let $\psi = (\pi \times \theta) \ind_{S_r \times S_s}^{S_{r+s}}$. 
It follows from Lemma~\ref{lemma:lr} that if $\chi^\nu$ is a
minimal constituent of $\psi$ then there are partitions $\lambda$
and $\mu$ such that $\left< \pi, \chi^\lambda \right> \ge 1$, 
$\left< \theta, \chi^\mu \right> \ge 1$ and $\nu = \lambda \cup \mu$.

Suppose that $\chi^\lambda$ is not a minimal constituent of $\pi$.
Then there exists
a partition $\lambda^\star$ such that $\lambda^\star \lhd \lambda$ and
$\left< \pi, \chi^{\lambda^\star} \right> \ge 1$. 
By Lemma~\ref{lemma:lr} we have 
\[ \left<\psi, \chi^{\lambda^\star \cup \mu} \right> \ge 1. \]
It is easily seen that $\lambda^\star \cup \mu \lhd \lambda \cup \mu$;
this contradicts the minimality of~$\chi^\nu$. Therefore $\chi^\lambda$
is a minimal constituent of $\pi$ and similarly, 
$\chi^\mu$ is a minimal constituent of $\theta$.
\end{proof}

We are now ready to prove Proposition~\ref{prop:gen}. Let 
$N = m(1)n(1) + \cdots + m(t)n(t)$.
Since the
$m(r)$ are distinct, 
\[ \phi^\mu = 
\Bigr(\phi^{(m(1)^{n(1)})} \times \cdots \times \phi^{(m(t)^{n(t)})}
\Bigr) \Ind_{S_{m(1)n(1)} \times \cdots \times S_{m(t)n(t)}}^{S_N}.
\]
The proposition now follows by repeated applications of Lemma~\ref{lemma:join}.

We finish with the observation that 
the converse to Proposition~\ref{prop:gen} (and to Lemma~\ref{lemma:join})
is false. This can be demonstrated using Corollary~\ref{cor:mins} in 
\S 5.2. It follows from this corollary that $\phi^{(5^5)}$ 
has 
%\begin{align*}
\[ \chi^{(5^5) \star (2,1,1)} = \chi^{(7,6,6,4,2)} \quad\text{and} \quad %\\
\chi^{(5^5) \star (1^4)} = \chi^{(6,6,6,6,1)}
\]
%\end{align*}
as minimal constituents. 
Similarly, $\phi^{(3^5)}$ has
%\begin{align*}
\[ \chi^{(3^5) \star (4)}= \chi^{(7,2,2,2,2)} \quad \text{and} \quad %\\
\chi^{(3^5) \star (3,1)} = \chi^{(6,4,2,2,1)} 
\]
%\end{align*}
as minimal constituents.
It is clear that
\[ (6,6,6,6,1) \cup (7,2,2,2,2) \rhd (7,6,6,4,2) \cup (6,4,2,2,1). \]
Hence $\chi^{(6,6,6,6,1) \cup (7,2,2,2,2)}$ 
is not a minimal constituent of $\phi^{(5^5,3^5)}$, even though
it arises from the $\cup$-construction applied to a minimal
constituent of $\phi^{(5^5)}$ and a minimal constituent of $\phi^{(3^5)}$.

\def\cprime{$'$} \def\Dbar{\leavevmode\lower.6ex\hbox to 0pt{\hskip-.23ex  \accent"16\hss}D}

\bigskip

{\footnotesize
\begin{tabular}{lr}
 Rowena Paget \\
 School of Mathematics, Statistics and Actuarial Science \\
 University of Kent
Canterbury  \\
 CT2 7NF, UK   \\
{\tt  R.E.Paget@kent.ac.uk} \\
\end{tabular}
}

\medskip
{\footnotesize
\begin{tabular}{lr}
Mark Wildon \\
Department of Mathematics  \\
Royal Holloway, University of London \\
Egham \\
TW20 0EX, UK \\
{\tt mark.wildon@rhul.ac.uk}
\end{tabular}
}


\begin{thebibliography}{1}
\bibitem{Brion}Michel Brion. \newblock Stable properties of plethysm: on two conjectures of Foulkes. \newblock {\em Manuscripta Math.}, 80:347--371, 1993.

\bibitem{CR}Charles.~W. Curtis and Irving Reiner. \newblock {\em Representation theory of finite groups and associative algebras}, 
\newblock reprint of the 1962 original, AMS Chelsea Publishing, Providence, RI, 2006.


\bibitem{DentSiemons}Suzie~C. Dent and Johannes Siemons. \newblock On a conjecture of {F}oulkes. \newblock {\em J. Algebra}, 226(1):236--249, 2000.

\bibitem{Foulkes}H.~O. Foulkes. \newblock Concomitants of the quintic and sextic up to degree four in the  coefficients of the ground form. \newblock {\em J. London Math. Soc.}, 25:205--209, 1950.

\bibitem{FultonYT}William Fulton. \newblock {\em Young tableaux}, volume~35 of {\em London {M}athematical  {S}ociety student texts}. \newblock CUP, 1997.

\bibitem{Howe}R. Howe. \newblock $(\mathrm{GL}_n, \mathrm{GL}_m)$-duality and symmetric plethysm. \newblock {\em Proc. Indian Acad. Sci. (Math. Sci.)},  97:85--109, 1987.

\bibitem{James}G.~D. James. \newblock {\em The representation theory of the symmetric groups}, volume 682  of {\em Lecture Notes in Mathematics}.
\newblock Springer, Berlin, 1978.

\bibitem{McKay1}Tom McKay. \newblock On plethysm conjectures of {S}tanley and {F}oulkes. \newblock {\em J. Algebra}, 319(5):2050--2071, 2008.

\bibitem{MN}Jurgen M{\"u}ller and Max Neunh{\"o}ffer. \newblock Some computations regarding {F}oulkes' conjecture.
\newblock {\em Experiment. Math.}, 14(3):277--283, 2005.

\bibitem{Stanley}R.~P.~Stanley. \newblock 
Positivity problems and conjectures in algebraic combinatorics, in: 
Mathematics: Frontiers and Perspectives,
295--319. \newblock Amer. Math. Soc., 2000. 

\bibitem{Thrall}R.~M. Thrall. \newblock On symmetrized {K}ronecker powers and the structure of the free {L}ie  ring. \newblock {\em American Journal of Mathematics}, 64(2):371--388, 1942.\end{thebibliography}
\end{document}